\documentclass{amsart}

\title[Perturbed Basins]{Perturbed Basins of Attraction}
\author{Han peters}
\date{November 15, 2004}

\usepackage{amscd}

\newtheorem{theorem}{Theorem}
\newtheorem{lemma}[theorem]{Lemma}
\newtheorem{proposition}[theorem]{Proposition}

\newtheorem{conjecture}[theorem]{Conjecture}

\theoremstyle{definition}
\newtheorem{defin}[theorem]{Definition}

\theoremstyle{remark}
\newtheorem{remark}[theorem]{Remark}

\newcommand{\NN}{\mathbb{N}}
\newcommand{\RR}{\mathbb{R}}
\newcommand{\CC}{\mathbb{C}}
\newcommand{\PP}{\mathbb{P}}
\newcommand{\aut}{\mathrm{Aut}_0(\CC^k)}

\def\B{{\bf B}}

\def\N{{\mathcal{N}}}
\begin{document}
\maketitle

\begin{abstract}
Let $F$ be an automorphism of $\CC^k$ which has an attracting
fixed point. It is well known that the basin of attraction is
biholomorphically equivalent to $\CC^k$. We will show that the
basin of attraction of a sequence of automorphisms $f_1, f_2,
\ldots$ is also biholomorphic to $\CC^k$ if every $f_n$ is a small
perturbation of the original map $F$.
\end{abstract}

\section{introduction}

 We are interested in the following conjecture, posed by Eric
 Bedford:

\begin{conjecture}\label{bedford}
 Let $F$ be a holomorphic automorphism of a complex manifold $X$,
 which is hyperbolic on an invariant compact subset $K$. Then for every $p
 \in K$, the stable manifold through $p$ is biholomorphically equivalent to
  complex Euclidean space.
\end{conjecture}

It was proved by Jonsson and Varolin \cite{jv} that Conjecture
\ref{bedford} holds almost everywhere with respect to any
invariant probability measure supported on $K$.

One could show that Conjecture \ref{bedford} holds by solving a
related problem in non-autonomous dynamics. For $z \in K$, let
$\phi_z$ be a biholomorphic mapping that maps the local stable
manifold through $z$ onto the local stable tangent space at $z$
(which after scaling we can identify with the unit ball in
$\CC^k$), where $k$ is the dimension of the stable manifold. We
may assume that $\phi_z^\prime(z)$ equals the identity. Let $p_0 =
p, p_1, p_2, \ldots$ be the orbit of $p \in K$. Then the maps
$$
f_n (z) = \phi_{p_{n+1}} F \phi_{p_n}^{-1},
$$
are biholomorphic mappings from the unit ball into the unit ball
(we may assume that these maps are really contracting, since
otherwise we consider some high iterate of $F$ instead). The maps
satisfy the conditions of the main theorem in \cite{fst}, and thus
there is a basin of attraction, defined in terms of the
{\emph{tail space}}. It is clear from the definition that this
basin of attraction is biholomorphic to the stable manifolds, and
the main theorem of \cite{fst} says that the basin of attraction
is biholomorphic to a basin of attraction of a sequence of global
automorphisms, which are uniformly contracting on some
neighborhood of the origin. Thus, Conjecture \ref{bedford} can be
answered positively by proving that the following conjecture
holds:

\begin{conjecture}\label{compact}
Let $f_1, f_2, \ldots$ be a sequence of automorphisms of $\CC^k$
that fix the origin. Suppose that there exist $a,b \in \RR$ that
satisfy $0 < a< b < 1$, and such that the following holds for
every $ n\in \NN$ and every $z \in \B$:
$$
a\|z\| \le \|f_n(z)\| \le b \|z\|.
$$
Then the basin of attraction of the sequence $f_1, f_2, \ldots$ is
biholomorphic to $\CC^k$.
\end{conjecture}

The {\emph{basin of attraction}} of a sequence $f_1, f_2 \ldots$
of automorphisms that all fix the origin is defined as the set of
all points $z \in \CC^k$ such that $f(n)(z) = f_n \circ \cdots
\circ f_1 (z) \rightarrow 0$.

If the automorphisms $f_1, f_2 \ldots$ do not all fix the origin,
but there does exist a small neighborhood $\N$ of $0$ such that
$f_n(\N) \subset \subset \N$ for every $n \in \NN$, then we define
the basin of attraction of the sequence $f_1, f_2, \ldots$ as
$\bigcup f(n)^{-1} \N$.

If all the maps $f_n$ are equal then the conjecture follows
immediately from the following theorem, which follows from the
work of Sternberg \cite{sg}, and was proved independently by Rosay
and Rudin \cite{rr}:

\begin{theorem}\label{RR} Let $F$ be an automorphism of
$\CC^k$, which has an attracting fixed point at $0$, then the
basin of attraction of $F$ is biholomorphically equivalent to
$\CC^k$.
\end{theorem}

Here, $F$ is attracting at the origin means that all eigenvalues
of $F^\prime(0)$ are smaller than $1$ in absolute value. We will
focus on Conjecture \ref{compact}. Our main result is the
following generalization of Theorem \ref{RR}:\\
\\
{\bf{Main Theorem.}} {\emph{ Let $F$ be an automorphism of $\CC^k$
which has an attracting fixed point at $0$. Then there exists an
$\epsilon>0$ such that for any sequence $f_1, f_2, \ldots$ of
automorphisms which satisfy $\|F(z), f_n(z)\| < \epsilon$ for any
$n \in \NN$ and all $z$ in the unit ball, one has that the basin
of attraction of the sequence $f_1, f_2, \ldots$ is
biholomorphic to $\CC^k$.}}\\

The maps $f_n$ don't necessarily have $0$ as an attracting fixed
point. However, if we choose $\epsilon$ small enough then there
exists a small neighborhood of $\N$ of $0$, such that $f_n(\N)
\subset \subset \N$ for all $n \in \NN$. Notice that the
definition of the basin of attraction is independent of the choice
of $\N$. Also notice that if the maps $f_n$ do fix the origin,
then the two definitions of the basin of attraction are exactly
equal.

Hence the Main Theorem gives sufficient conditions for when the
basin of attraction of a sequence $f_1, f_2, \ldots$ of
automorphisms which have an attracting fixed point at the origin
is biholomorphic to $\CC^k$. The following theorem can be found in
\cite{wd}:
\begin{theorem}\label{square}
Let $f_1, f_2, \ldots$ be a sequence of automorphisms of $\CC^k$
which all fix the origin, and let $a,b \in \RR$ that satisfy $0 <
a < b < 1$ and the additional condition $b^2 < a$. Suppose that
the following holds for every $ n\in \NN$ and every $z \in \B$:
$$
a\|z\| \le \|f_n(z)\| \le b \|z\|.
$$
Then the basin of attraction of the sequence $f_1, f_2, \ldots$ is
biholomorphic to $\CC^k$.
\end{theorem}

The proof of Theorem \ref{square} can be adapted to get the
following result which is slightly stronger:

\begin{theorem}\label{blocks}
Let $f_1, f_2, \ldots$ be a sequence of automorphisms which all
fix the origin, and let $t \in (1,2)$. Suppose that for every $
n\in \NN$ there exists $c_1, c_2, \ldots \in (0,1)$ such that
$\prod c_n = 0$ and the following holds for every $z \in \B$:
$$
c_n^t\|z\| \le \|f_n(z)\| \le c_n \|z\|.
$$
Then the basin of attraction of the sequence $f_1, f_2, \ldots$ is
biholomorphic to $\CC^k$.
\end{theorem}

The idea in the proof of both Theorem \ref{square} and
\ref{blocks} is the following:  one defines the maps $H_n :=
A(n)^{-1}f(n)$, where $A(n) = A_n \circ \cdots \circ A_1$ and $A_n
= f_n^\prime(0)$. As $n \rightarrow \infty$, one gets that $H_n$
converges uniformly on compact subsets of the basin of attraction
to a biholomorphic mapping from the basin onto $\CC^k$.

The following theorem was proved in \cite{pw} and is also a
generalization of Theorem \cite{rr}:

\begin{theorem}\label{repeat}
Let $f_1, f_2, \ldots$ be a sequence of automorphisms of $\CC^k$,
which are all attracting at the origin. Then there exist large
enough integers $n_1, n_2, \ldots$ such that the basin of
attraction of the maps $f_1^{n_1}, f_2^{n_2} \ldots$ is
biholomorphic to $\CC^k$.
\end{theorem}

Here the definition of the basin of attraction must be slightly
changed for this theorem to hold. The exact definition is given in
\cite{pw}, as well as an example that shows that Theorem
\ref{repeat} does not hold if one uses a more straightforward
definition.

In the next section we will set our notation, and introduce the
main ingredient in the proof of our Main Theorem, namely Theorem
\ref{stable}. In the Third section we will prove our Main Theorem
using Theorem \ref{stable}. In section 4 we give the main idea of
the proof of Theorem \ref{stable}, and we finish up by proving
Theorem \ref{stable} in the fifth section.

\section{Preliminaries}

Throughout the paper we write $\aut$ for the set of automorphisms
of $\CC^k$ that fix the origin. We denote by $B(r)$ the ball or
radius $r$ centered at the origin and we write $\B$ for the unit
ball. For a sequence $F_1, F_2, \ldots$ we will use the notation
$F(n)$ for the composition of the first $n$ maps, $F_n \cdots F_1$
and we write $F(m,n) = F_n \cdots F_m = F(n)F(m)^{-1}$.

It follows from a well known fact from linear algebra (called $QR$
factorization) that for any automorphism $F \in \aut$ we can find
a unitary matrix $U \in U_k(\CC)$ such that $(U \circ F)^\prime
(0)$ is lower triangular. The map $U$ is in fact unique up to
multiplying rows with complex numbers of unit length. We say that
$U F$ is in {\emph{lower triangular form}}.

\begin{remark}\label{QR}
Given some sequence $F_1, F_2, \ldots \in \aut$ we can find
unitary matrices $U_1, U_2 \ldots$ such that $U_1 F_1, U_2 F_2
U_1^{-1}, U_3 F_3 U_2^{-1}, \ldots$ are all in lower triangular
form.

If all the maps $F_n$ are attracting at the origin, then we define
$\tilde{F}_n = U_n F_n U_{n-1}^{-1}$. We have that $\tilde{F}(n) =
U_n F(n)$, therefore the basin of attraction of the sequence
$\{\tilde{F}_n\}$ is exactly equal to the basin of the original
maps, so we may assume that our original maps were already in
lower triangular form.
\end{remark}

We will say that $F \in \aut$ which is in lower triangular form is
\emph{correctly ordered} if $F^\prime(0)$ has diagonal entries
(from upper left to lower right) $\lambda_1, \ldots , \lambda_k$
that satisfy the following condition:
$$
|\lambda_j| |\lambda_i| < |\lambda_l|,
$$
for $l \le j$ and any $i$. Note in particular that $|\lambda_j| <
1$ for every $j$, and that $|\lambda_j|^2 < |\lambda_l|$ for $ j
\ge l$. Note also that if the eigenvalues satisfy the ordering
$|\lambda_1| \ge \cdots \ge |\lambda_k|$ then $F$ is correctly
ordered, but the sequence $|\lambda_1|, \ldots ,|\lambda_k|$ may
fail to be decreasing by a relatively small amount.

\begin{defin}\label{attracting}
Let $\mathcal{F} \subset \aut$ be a family of correctly ordered
automorphisms. We say that $\mathcal{F}$ is {\emph{uniformly
attracting}} if there exist $a, b \in \RR$ with $0 < a < b < 1$
and we have that for every $F \in \mathcal{F}$ and every $z \in
\B$,
$$
a \|z\| \le \|F(z)\| \le b \|z\|.
$$
Additionally we require that there exists a uniform $\xi <1$ such
that
\begin{eqnarray} \label{eigenvalues}
|\lambda_i||\lambda_j| \le \xi |\lambda_l|,
\end{eqnarray}
for $l \le j$ and any $i$.
\end{defin}

We will prove the following theorem:
\begin{theorem} \label{stable}
Let $F_1, F_2, \ldots \in \aut$ be a uniformly attracting sequence
of correctly ordered automorphisms. Then the basin of attraction
of $F_1, F_2, \ldots$ is biholomorphic to $\CC^k$.
\end{theorem}

\section{Proof of the Main Theorem}

We will now prove our main Theorem, using Theorem \ref{stable}.
The proof consists of several (but finitely many ) steps, and in
each step we may decrease the value of $\epsilon$.

We may assume that $F$ is such that $F^\prime(0)$ is a lower
triangular matrix with diagonal entries (from upper left to lower
right) $\lambda_1, \ldots ,\lambda_k$ that satisfy $|\lambda_1|
\ge \ldots \ge |\lambda_k|$. We may also assume that the
off-diagonal terms of $F^\prime(0)$ are arbitrarily small.

Now we show that we may assume that all the maps $f_n$ have an
attracting fixed point at $0$. Since $F^\prime(0)$ is arbitrarily
close to a diagonal matrix, we choose the neighborhood of the
origin $\N$ so small that $F$ is contracting on $\N$, i.e. that
there exist some $\theta < 1$ such that for all $x,y \in \N$ we
have that $\|F(x) - F(y)\| < \theta \|x - y\|$. Then we can make
$\epsilon$ small enough such that every $f_n$ is also contracting
on $\N$ (with uniform constant $\theta < 1$ that is slightly
increased if necessary). Let $x_0 = 0, x_1, x_2, \ldots$ be the
orbit of $0$, i.e. $x_n = f_n(x_{n-1})$. We have that
$f_n(\mathcal{N}) \subset \subset \mathcal{N}$ for every $n \in
\NN$, so that $x_n \in \mathcal{N}$ for every $n \in \NN$. Let
$T_n$ be the translation of $\CC^k$ that maps $x_n$ to $0$. Then
define $\tilde{f_n} = T_n \circ f_n \circ T_{n-1}^{-1}$. We have
that $\tilde{f_n}(0) = 0$ for all $n$, $0$ is an attracting fixed
point for every map $f_n$, and the maps $\tilde{f_n}$ are still
arbitrarily close to the original map $F$. Since $\tilde{f}(n) =
T_n f(n)$, we have that the approximate basin of attraction of the
sequence $f_1, f_2, \ldots$ is exactly equal to the basin of
attraction of $\tilde{f}_1, \tilde{f}_2, \ldots$. Therefore, we
may as well assume that all the maps $f_n$ have $0$ has an
attracting fixed point, and we only consider the basin of
attraction of the sequence $f_1, f_2, \ldots$.

We will now prove the Main Theorem in the case that $k = 2$. We
may assume that $|\lambda_2|$ is strictly smaller than
$|\lambda_1|$, otherwise the result follows easily from Theorem
\ref{square}. The fact that $F^\prime(0)$ is lower diagonal means
exactly that $(0,1)$ is an eigenvector of $F^\prime(0)$. Let
$\Phi$ be the action on $\PP^1$ induced by the mapping
$F^\prime(0)$. It follows that $\Phi([0:1]) = [0:1]$, and the
multiplier of $\Phi$ at $[0:1]$ is exactly
$\frac{\lambda_1}{\lambda_2}$. Hence $[0:1]$ is a repelling fixed
point, and there exists an arbitrarily small neighborhood of
$[0:1]$, say $\mathcal{N}$ such that $\N \subset\subset \Phi(\N)$.

Let $\phi_n$ be the action on $\PP^1$ induced by $f_n^\prime(0)$.
We can make sure that $\phi_n$ is arbitrarily close to $\Phi$ such
that $\N \subset\subset \phi_n(\N)$ for all $n \in \NN$. Then we
have that $\N \supset \supset \phi_1^{-1}(\N) \supset\supset
\phi(2)^{-1}(\N) \cdots$. It follows that
$$
\bigcap_{n \in \NN} \phi(n)^{-1}(\N) \neq \emptyset.
$$

Let $v \in \bigcap \phi(n)^{-1}(\N)$, and let $v_0 = v, v_1, v_2,
\ldots$ be the orbit of $v$, i.e. $\phi_n(v_{n-1})= v_n$. We have
that $v_n \in \N$ for every $n \in \NN$.

Let $U_n$ be a unitary $2 \times 2$ matrix that maps some length
$1$-representative of $v_n$ in $\CC^2$ onto $(0,1)$, and define
$\tilde{f}_n := U_n f_n U_{n-1}^{-1}$. Then it follows that
$(0,1)$ is an eigenvector of every map $\tilde{f}_n^\prime(0)$,
and thus we have that $\tilde{f}_n^\prime(0)$ is lower triangular.
The basin of the sequence $\tilde{f}_1, \tilde{f}_2, \ldots$ is
exactly equal to $U_0(\Omega)$, so in particular biholomorphic to
$\Omega$. Since the unitary matrices $U_n$ are arbitrarily close
to the identity matrix, we have that $\tilde{f}_n$ is arbitrarily
close to $F$, and it follows that the sequence $\tilde{f}_1,
\tilde{f}_2, \ldots$ satisfy the properties in Theorem
\ref{stable} and we are done.\\

We will use a similar argument to that for $2$ dimensions to prove
the Main Theorem in the general case.

\begin{lemma}\label{grass}
Let $A$ be a lower triangular $k \times k$ matrix whose diagonal
entries $\lambda_1, \lambda_2, \ldots \lambda_k$ satisfy
$|\lambda_1| \ge \cdots \ge |\lambda_k|$. Suppose that $l \in \{1,
\ldots k-1\} $ is such that $|\lambda_l| > |\lambda_{l+1}|$, and
let $L$ be the linear subspace of $\PP^{k-1}$ defined by
$$
L = \{ [z_1:\cdots : z_k] \in \PP^{k-1} \mid z_1 = \cdots = z_{l}
= 0\}.
$$
Then there exists an arbitrarily small neighborhood $\N \subset
G(k-l-1, k-1)$ of $L$ such that $A(\N) \supset \supset \N$.
\end{lemma}
\begin{proof}
First assume that $A$ is diagonal. For $X \in G(k-l-1, k-1)$ close
enough to $L$ we can write
$$
X = \{ [z_1:\cdots : z_k] \in \PP^{k-1} \mid z_1 = \sum_{j \ge l}
\epsilon_{1,j} z_j, \ldots ,z_l = \sum \epsilon_{l,j} z_j\}.
$$
We define $\mathrm{d}(X, L) = \max |\epsilon_{i,j}|$.

For $\delta >0$ small let $\N_{\delta} = \{X \in G(k-l-1, k-1)
\mid \mathrm{d}(X,L) < \delta$\}. Since $A$ is diagonal and
$|\lambda_l| > |\lambda_{l+1}|$ for all $i \le l$ and $j \ge l+1$,
we have that
$$
A(\N_{\delta}) \supset
\N_{\frac{|\lambda_l|}{|\lambda_{l+1}|}\delta} \supset \supset
\N_{\delta}
$$
Fix some $\delta$. Then $\N_{\delta}$ will also work for
arbitrarily small perturbations of the diagonal matrix $A$.

We can conjugate any $A$ with an invertible linear mapping $T$
such that $T^{-1} A T$ is lower triangular and arbitrarily close
to a diagonal mapping (whose diagonal entries are exactly those of
$A$). Hence the set $T(\N_\delta)$ will suffice for $A$.
\end{proof}

Let $l \in \{1, \ldots ,k-1\}$ be such that $|\lambda_{l}| >
|\lambda_{l+1}| = |\lambda_k|$, and let $\Phi$ be the action of
$F^\prime(0)$ on $\PP^{k-1}$ (as in the 2-dimensional case, we may
assume that such an $l$ exists, otherwise we are done by Theorem
\ref{square}). Then it follows from Lemma \ref{grass} that there
exists a small neighborhood $\N \in G(k-l-1, k-1)$ of $L$ (where
$L$ is as in Lemma \ref{grass}, such that $\Phi(\N) \supset\supset
\N$. As in the 2-dimensional case, we denote by $\phi_n$ the
action of $f_n$ on $\PP^{k-1}$. Then we can decrease the value of
$\epsilon$ if necessary to get that $\phi_n(\N) \supset\supset \N$
for every $n \in \NN$. Therefore we can find an orbit of linear
subspaces $L_0, L_1, \ldots$ in $\N$ with $\phi_n(L_{n-1}) = L_n$
for every $n \in \NN$.

We denote by $U_n$ a unitary matrix arbitrarily close to the
identity that maps a representative of $L_n$ in the unit sphere in
$\CC^k$ onto the set
$$
T = \{ z \in \CC^k \mid \|z\| = 1, z_1 = \cdots = z_l = 0 \}.
$$
We define the maps $\tilde{f}_n$ by
$$
\tilde{f}_n = U_n f_n U_{n-1}^{-1}.
$$
Then we have that for every $v \in T$ there exists a $c \in (0,1)$
such that $\tilde{f}_n^\prime(0)(v) \in c T$, and therefore we
have that for any $n \in \NN$ the entries in the $p$-th rows and
$q$-th columns of the matrix $f_n^\prime(0)$ are equal to $0$ for
$p \le l$ and $q \ge l+1$. Since the matrices $U_n$ are
arbitrarily close to the identity, we have that
$\tilde{f}_n^\prime(0)$ is arbitrarily close to $F^\prime(0)$.

 We may also assume that the $(k-l) \times (k-l)$ blocks
  in the lower right corner of the matrices
$\tilde{f}_n^\prime(0)$ are lower
 triangular, since we can apply QR-factorizations as in Remark \ref{QR}
  to these $(k-l) \times (k-l)$
blocks. The diagonal entries in these blocks must in absolute
value be arbitrarily close to $|\lambda_k|$, since in this last
step have only composed with unitary matrices, and in Step 2 we
made sure that $F^\prime(0)$ is arbitrarily close to a diagonal
matrix. Also, it follows that the off-diagonal terms in the $(k-l)
\times (k-l)$ block in the lower right corner of the matrices
$\tilde{f}_n^\prime(0)$ are arbitrarily small.

 Since the last $k-l$
columns of the matrices $\tilde{f}_n^\prime(0)$ are already lower
triangular, and the corresponding eigenvalues are the smallest in
absolute value, we can restrict ourselves to the first $l$
dimensions and apply the same arguments to the next eigenvalues of
$F^\prime$ which are equal in absolute value.

In finitely many steps we get that all the maps
$\tilde{f}_n^\prime(0)$ are lower triangular matrices. To change
from $f$ to $\tilde{f}$, we only compose with unitary matrices
that are arbitrarily close to the identity. As a result, we get
that the diagonal entries of the maps $\tilde{f}_n^\prime(0)$ are
arbitrarily close to the diagonal entrees of $F^\prime(0)$.
Therefore, all the maps in the sequence $f_1, f_2, \ldots$ are
correctly ordered.

We have constructed unitary matrices $U_0, U_1, \ldots$ such that
for every $n \in \NN$ we have that $\tilde{f}(n) = U_n f(n)
U_0^{-1}$. It follows that the basin of attraction of the sequence
$\tilde{f}_1, \tilde{f}_2, \ldots$ is equal to the image of the
basin of attraction of the original sequence $f_1, f_2 \ldots$
under the map $U_0^{-1}$, which is a biholomorphism. The sequence
$\tilde{f}_1, \tilde{f}_2, \ldots$ satisfies the conditions in
Theorem \ref{stable}, and the Main Theorem.

\section{Main Ideas of Theorem \ref{stable}}

In the proof of Theorem \ref{stable} we will construct a sequence
of maps $\Psi_n := G(n)^{-1} \circ X_n \circ F(n)$ which converges
to the Fatou Bieberbach mapping $\Psi : \Omega \rightarrow \CC^k$.
Here the $G_n$'s are \emph{lower triangular polynomial mappings}
(as in \cite{rr}), and the $X_n$'s are polynomial mappings whose
linear parts are the identity. To be more specific, we will start
with some choice for $X_1$ and then define $X_n = [G_n \circ
X_{n-1} \circ F_n^{-1}]_{p}$. Here $d$ is some large integer and
we mean by $[\cdot]_d$ that we discard all terms of degree $d+1$
and higher. It follows immediately that

\begin{eqnarray} \label{difference}
\|G_n^{-1}X_nF_n(z) - X_{n-1}(z)\| = O(\|z\|^{d+1}).
\end{eqnarray}

The challenge is to choose a bounded sequence of lower triangular
polynomial mappings $\{G_n\}$ such that there exists a bounded
orbit $\{X_n\}$, and to get the integer $d$ as large as necessary.
We first show that we can do this in some simpler cases before we
complete the proof of theorem \ref{stable}.

Recall that a polynomial selfmap $G= (g_1, \ldots g_k)$ of $\CC^k$
with $G(0) = 0$ is called lower triangular if
$$
g_j(z) = c_j z_j + h_j (z_1, \ldots z_{j-1}),
$$
for all $j \in \{1, \ldots , k\}$.

\begin{lemma} \label{upper}
Let $G_1, G_2, \ldots$ be a sequence of lower triangular
polynomial mappings of some fixed degree, whose coefficients are
uniformly bounded. Then we have the following:

(a) The degrees of the maps $G(n)$ are bounded, and there is a
constant $\beta < \infty$ so that
$$
G(n)(\B) \subset B(\beta^n)
$$

(b) If also $|c_i|< \theta < 1$ for all $G_j$, then $G(n)(z)
\rightarrow 0$, uniformly on compact subsets of $\CC^k$, and for
every $R > 0$ and $\epsilon
>0$ there exist an $N \in \NN$ such that for $n \ge N$ we have
that
$$
G(n)(B(R)) \subset B(\epsilon).
$$
\end{lemma}
This lemma is similar to Lemma 1 in the appendix of \cite{rr}, and
so is its proof.

We now prove the following simple lemma.
\begin{lemma} \label{affine}
Let $F_1, F_2, \ldots$ be a uniformly bounded sequence of affine
maps on $\CC$ that are all expanding. Then there exists a bounded
orbit $z_0, z_1, z_2, \ldots$
\end{lemma}

Let $F_n (z) = a z + b$. With uniformly bounded we mean that the
constants $a$ and $b$ are bounded from above, and that the
constants $a$ are also bounded from below by some $c > 1$.

\begin{proof}
Since the sequence is compact, there exists a constant $R$ such
that
\begin{eqnarray*}
F_i(B(R)) \supset \supset B(R)
\end{eqnarray*}
for all $i \in \NN$. Now we have that
\begin{eqnarray*}
B(R) \supset\supset F(1)^{-1}B(R) \supset\supset F(2)^{-1}B(R)
\supset\supset \ldots,
\end{eqnarray*}
So the intersection $\bigcap F(n) ^{-1} \Delta (R)$ contains a
point $z_0$. Then we have that the sequence $z_0, z_1, z_2,
\ldots$, where $z_{n} = F_n(z_{n-1})$, is a bounded orbit. Indeed,
$|z_n| < R$ for all $n \in \NN$.
\end{proof}

\subsection{Bounded Orbits in $\CC$}

The argument that shows that there exist bounded sequences $G_n$
and $X_n$ is somewhat complicated. To make it more clear we will
first prove the existence in the one dimensional case, in which
case it is much easier.

Let $f_1, f_2, \ldots$ be a sequence of polynomials with uniformly
bounded coefficients of the form $f_n (z) = \lambda_n z + h.o.t.$,
where $|\lambda_n|
> \theta
> 1$. Let $g_n(z) = \lambda_n^{-1} z$. Let $\mathcal{P}_d$ be
the space of polynomials of the form $z+ h.o.t.$ of degree at most
$d$. Define the map $\phi_n:\mathcal{P}_d \rightarrow
\mathcal{P}_d$ by $\phi_n(X) = [g_n\circ X \circ f_n]_d$, where we
mean by $[\cdot ]_d$ that we discard all terms of degree strictly
higher than $d$.

\begin{proposition}
There exists a bounded sequence $X_0, X_1, \ldots \in
\mathcal{P}_d$ with $X_n = \phi_n(X_{n-1})$.
\end{proposition}
\begin{proof}

We will use induction on the degree. For degree $1$ the statement
is clear since we take $X_n (z) = z$ for all $n$.

Suppose we have constructed a bounded orbit $Y_n$ that satisfies
condition \eqref{difference} for some $p = d-1 \ge 1$. We will add
terms of degree $p$ so that the sequence satisfies condition
\eqref{difference} for $p = d$. Say $X_0 = Y_0 + c_0 z^p$ for some
$c_0 \in \CC$. Then it is easy to see that $X_1 = \phi_1 (X_0) =
Y_1 + c_1 z^d$, where $c_1$ is equal to $\lambda_1^{d-1} c_0$ plus
some constant depending linearly on the coefficients of $Y_0$ and
the higher order terms of $f_1$ (which are uniformly bounded by
assumption). We see that we get a uniformly bounded sequence of
expanding affine maps $\tilde\psi_n : \CC \rightarrow \CC$ that
take $c_{n-1}$ to $c_n$. It now follows from Lemma \ref{affine}
that we can choose $c_0$ such that the sequence $c_0, c_1, \ldots$
is bounded.

It follows by induction that we can find a bounded sequence $X_n$
for any $d \in \NN$.
\end{proof}

\subsection{Degree $2$ polynomials in $\CC^2$}

Next we show that we can get bounded sequences $\{G_n\}$ and
$\{X_n\}$ in the two dimensional case for degree $2$, where we can
explicitly calculate the maps. We write $F_n (x, y) = (\lambda_n x
, \mu_n y + a_n x) + h.o.t.$ and $(F_n^{-1}) (x, y) =
(\lambda_n^{-1} x , \mu_n^{-1} y + b_n x) + h.o.t.$ and let
$G_n(x,y) = (\lambda_n x , \mu_n y + a_n x + d_n x^2)$ for some
constants $d_n \in \CC$ to be chosen later. We will also write
$X_n(x,y) = (x + \alpha_n y^2 + \beta_n xy + \gamma_n x^2, y +
\delta_n y^2 + \epsilon_n xy + \zeta_n x^2)$. We will identify the
map $X_n$ with $(\alpha_n, \beta_n, \gamma_n, \delta_n,
\epsilon_n, \zeta_n) \in \CC^6$. Consider the map $X_{n-1} \mapsto
X_n = [G_n \circ X_{n-1} \circ F_n^{-1}]_2$. We have that
\begin{align*}
\alpha_n &= \lambda_n \mu_n^{-2} \alpha_{n-1} + l_{1,n}, \\
\beta_n &= \mu_n^{-1} \beta_{n-1} + l_{2,n}(\alpha_{n-1}),\\
\gamma_n &= \lambda_n^{-1} \gamma_{n-1} + l_{3,n}(\alpha_{n-1},
\beta_{n-1}),\\
\delta_n &= \mu_n^{-1} \delta_{n-1} + l_{4,n}(\alpha_{n-1}),\\
\epsilon_n &= \lambda_n^{-1} \epsilon_{n-1} + l_{5,n}(\alpha_{n-1}, \beta_{n-1}, \delta_{n-1}),\\
\zeta_n &= \mu_n \lambda_n^{-2} \zeta_{n-1} +
l_{6,n}(\alpha_{n-1}, \ldots ,\epsilon_{n-1}) + d_n
\lambda_n^{-2}.
\end{align*}
Here the $l_{i,n}$ are linear maps that depend only on the
coefficients of $F_n^{-1}$ (which are uniformly bounded) and on
the given variables.

It follows from equation \eqref{eigenvalues} that $|\lambda_n
\mu_n^{-2}|
>1$ for any $n$, and therefore we get a uniformly bounded sequence of
expanding affine maps $\alpha_{n-1} \mapsto \alpha_n$. Hence it
follows from Lemma \ref{affine} that we can find $\alpha_0$ such
that the sequence $\alpha_0, \alpha_1, \ldots$ is bounded. Having
fixed the $\alpha_n$'s, we can use the same argument for the
$\beta_n$'s, since we also have that $|\mu_n^{-1}|>1$ for all $n$.
After we fix the $\beta_n$'s, we can find a bound on the
$\gamma_n$'s, then the $\delta_n$'s, and finally the
$\epsilon_n$'s.

We can't use the same argument for the $\zeta_n$'s, since we may
not have that $|\mu_n \lambda_n^{-2}| >1$. However, we can choose
$\zeta_0 = 0$, and then choose the constants $d_n$ such that
$\zeta_n = 0$ for every $n$. Hence we get the bounded sequences
$\{G_n\}$ and $\{X_n\}$.

\section{Proof of Theorem \ref{stable}}

The argument that we use to construct bounded sequences $\{X_n\}$
and $\{G_n\}$ for higher dimensions and higher degrees is
essentially the same as the argument for degree $2$ polynomial
mappings in $\CC^2$.

We will write $\lambda_{n,1}, \ldots , \lambda_{n,k}$ for the
diagonal terms of $F_n^\prime(0)$.

\begin{proposition}\label{bounded}
For any $d \ge 2$ we can find a bounded sequence of polynomial
mappings $X_0, X_1, \ldots$, where $X_n = I_k + h.o.t.$, and a
bounded sequence of lower triangular polynomial mappings $G_1,
G_2, \ldots$, where $G_n^\prime(0) = F_n^\prime(0)$, such that
\eqref{difference} holds for any $n \ge 1$.
\end{proposition}
\begin{proof}
We will construct bounded sequences $X_1, X_2, \ldots$ and $G_1,
G_2, \ldots$ such that the following equation holds:
\begin{eqnarray}\label{next}
X_n = [G_n X_{n-1} F_n^{-1}]_d,
\end{eqnarray}
for every $n \in \NN$. Notice that this will imply that Equation
\ref{difference} holds for every $n \in \NN$.

We write
$$
X_n(z) = (x_{n,1}(z), x_{n,2}(z), \ldots, x_{n,k}(z)),
$$
and
$$
x_{n,j} = \sum_{\alpha} c_{n,j,\alpha}z^\alpha,
$$
where $\alpha$ is a $k$-tuple, and $z^\alpha = z_1^{\alpha_1}
\cdots z_k^{\alpha_k}$. We refer to $c_{n,j,\alpha}z^\alpha$ as a
term of degree $|\alpha| = \alpha_1 + \cdots + \alpha_k$, index
$j$, and power $\alpha$. We similarly write $g_{n,j,
\alpha}z^{\alpha}$ for the term of the mapping $G_n$ of degree
$|\alpha|$, index $j$ and power $\alpha$.

For two $k$-tuples $\alpha$ and $\beta$ with $|\alpha| = |\beta|$
we will write $\alpha > \beta$ if $\alpha$ has a higher
lexicographical ordering than $\beta$. That is, $\alpha > \beta$
if and only if there is some $j \in \{1, \ldots , k\}$ such that
$\alpha_j < \beta_j$ and $\alpha_i = \beta_i$ for $i \in \{1,
\ldots j-1\}$ (so that $z_1^d$ comes first in the alphabet and
$z_k^d$ comes last).

We will use induction on the degree and index and reverse
induction on the power to fix all the terms of the sequences
$\{X_n\}$ and $\{G_n\}$.

Let $\alpha$ be some $k$-tuple, and let $j \in \{1, \ldots ,k\}$.
Suppose that we have fixed all the terms of degree up to
$|\alpha|-1$ in the sequences $\{X_n\}$ and $\{G_n\}$, such that
\eqref{next} holds for $d = |\alpha|-1$, and assume that the
corresponding coefficients are uniformly bounded. Also assume that
we have fixed all terms of degree $|\alpha|$ and index up to $j-1$
in the sequences $\{X_n\}$ and $\{G_n\}$ such that \eqref{next}
holds for index $1, \ldots , j-1$ and $d=|\alpha|$. Finally,
assume that all terms of degree $|\alpha|$, index $j$ and powers
$\beta > \alpha$ in the sequences $\{X_n\}$ and $\{G_n\}$ are
fixed such that \eqref{next} holds for the terms of degree
$|\alpha|$, index $j$ and power $\beta$.

We will continue to choose the terms of index $j$ and power
$\alpha$ in the sequences $\{X_n\}$ and $\{G_n\}$. It is clear
that after fixing the $c_{n, j, \alpha}$'s and $g_{n,j,\alpha}$'s,
Equation \eqref{next} will still hold for $d = |\alpha|-1$. Since
the linear parts of the $G_n$'s are lower triangular, it also
follows that \eqref{next} will still hold for degree $d =
|\alpha|$ and index $1, \ldots, j-1$. Furthermore, it follows from
the fact that $F_n^\prime(0)$ is lower triangular that
\eqref{next} will still hold for the terms of degree $d =
|\alpha|$, index $j$ and all powers $\beta$ with $\beta > \alpha$.

Once some $c_{n-1,j, \alpha}$ is chosen, it follows from
\eqref{next} that we must have
\begin{eqnarray}\label{constant}
c_{n, j, \alpha} = \lambda_{n,j} \lambda_n^{-\alpha} c_{n-1,j,
\alpha} + g_{n, j, \alpha} \lambda_n^{-\alpha}+ C_{n,j,\alpha}
\end{eqnarray}
where the constant $C_{n,j,\alpha}$ depends on the terms of
$f_{n}^{-1}$ (which are uniformly bounded for any fixed degree)
and the terms of $G_{n}$ and $X_{n-1}$ that we have already fixed
(which are also uniformly bounded per hypothesis).

If the term of index $j$ and power $\alpha$ is lower triangular
(i.e. if $\alpha_i = 0$ for $i \in \{j, \ldots k\}$), then we can
choose $c_{n,j, \alpha} = 0$ for all $n \in \NN$ and we choose
$g_{n, j, \alpha} = -C_{n, j, \alpha}\lambda_n^\alpha$ such that
\eqref{constant} holds for every $n \in \NN$. It follows that the
constants $g_{n,j, \alpha}$'s are bounded.

If the term of index $j$ and power $\alpha$ is not lower
triangular , then it follows from the hypothesis that the sequence
$\{F_n\}$ is uniformly attracting that $|\lambda_{n,j}
\lambda_n^{-\alpha}| > \xi > 1$ (where $\xi$ is as in Definition
\ref{attracting}). We also have that $g_{n,j, \alpha} = 0$ for all
$n \in \NN$, and it follows that the sequence of maps given by
$c_{n-1,j, \alpha} \mapsto c_{n,j, \alpha}$ is a uniformly bounded
sequence of affine maps on $\CC$ that are all expanding. It
follows from Lemma \ref{affine} that there exists a bounded
sequence $c_{n,j,\alpha}$.

So whether the term is lower triangular or not, we can always
choose bounded sequences $\{x_{n,j,\alpha}\}$ and
$\{g_{n,j,\alpha}\}$ such that \eqref{next} is satisfied.

The Proposition follows by induction.
\end{proof}

Now let $p$ be so large that for every $n \in \NN$ the eigenvalues
$\lambda_{n,1}, \ldots \lambda_{n,k}$ of $F_n^\prime(0)$ satisfy
$\xi |\lambda_{n,i}|^p < |\lambda_{n,j}|$ for $i,j \in \{1, \ldots
, k\}$, where $\xi > 1$. We can do this because the sequence
$\{F_n\}$ is uniformly attracting. We construct sequences
$\{G_n\}$ and $\{X_n\}$ as in Proposition \ref{bounded} for $d =
p$. It follows from part (a) of Lemma \ref{upper} that there
exists a $\gamma>1$ such that
$$
\|G(n)^{-1}(w) - G(n)^{-1}(w^\prime)\| \le \gamma^n \|w -
w^\prime\|,
$$
for any $w, w^\prime \in \B$. Recall that we assumed that there is
a constant $b < 1$ such that $\|F_n(z)\| < b \|z\|$ for any $z \in
\B$. Now fix an integer $q$ such that $\gamma b^q<\alpha<1$. We
change the mappings $X_n$ by adding higher order terms such that
\eqref{difference} holds for $d=q+1$. Since we chose $p \in \NN$
such that $|\lambda_i|^p < \xi |\lambda_j|$ holds for all $i,j$,
we have that $|\lambda_{n,j} \lambda_n^{-\alpha}| > \xi$ in
\eqref{constant} even for the terms that are lower triangular.
Hence we can make sure that the altered sequence $\{X_n\}$ is
bounded without changing the sequence $\{G_n\}$. Therefore we have
that
\begin{align}\label{done}
\|G_n^{-1}X_nF_n(z) - X_{n-1}(z)\| \le C\|z\|^{q+1},
\end{align}
for some $C > 0$ independent of $n \in \NN$ and every $z \in \B$.

The proof of Theorem \ref{stable} now follows quickly from an
argument similar to that of the proof of Theorem \ref{RR} which
can be found in the appendix of \cite{rr}.

Define the maps $\Psi_n : \Omega \rightarrow \CC^k$ by
\begin{eqnarray*}
\Psi_n := G(n)^{-1} \circ X_n \circ F(n).
\end{eqnarray*}
We will show that the maps $\Psi_n$ converge uniformly on compact
subsets of $\Omega$ to a biholomorphic map from $\Omega$ onto
$\CC^k$.

Since the sequence $\{X_n\}$ is bounded there is a radius $r$ such
that all the maps $X_n$ are invertible on $B(r)$.

Let $K$ be a compact set in $\Omega$. Then there is an $m \in \NN$
such that $F(m)(K) \subset B(r)$. Let $n \ge m$. Then we have that
$\|F(m,n)(z)\| \le b^{n-m} \|z\|$, for all $z \in B(r)$.

We notice that
\begin{align}\label{small}
\|\Psi_{n+1}(z) - \Psi_n(z)\| &= \| \left(G(n)^{-1}
G_{n+1}^{-1} X_{n+1} F_{n+1} - G(n)^{-1} X_n \right) F(n)(z)\|\\
&\le  \label{smaller} C \gamma^n (b^{n-m}\|F(m)(z)\| )^{q+1} <
\tilde{C} \alpha^n\|F(m)(z)\|^{q+1}.
\end{align}

Since the sequence $\{\alpha^n\}$ is summable, the maps $\Psi_n$
converge uniformly on compact subsets of $\Omega$ to a holomorphic
map $\Psi$. Also, for any compact subset $K$ of $\Omega$, there is
a large $N \in \NN$ such that for $n \ge N$, $\Psi_n$ is
biholomorphic on K. It is a well known fact that the limit of a
convergent sequence of biholomorphic mappings is either injective
or degenerate everywhere. We have that $\Psi_n^\prime (0) = I$ for
all $n$, therefore we have that  $\Psi^\prime(0) = I$ and the
limit map $\Psi$ is injective.

To prove surjectivity of the map $\Psi$, we may assume that we
have chosen $r$ so small that
$$
r^{q} \sum_{n=1}^\infty (C \alpha^n) < \frac{1}{2}.
$$
It follows from estimates \eqref{small} and \eqref{smaller} that
for $z \in B(r)$ and for any $n > m \ge 1$ we have that
$$
\|G(m,n)^{-1}X_n F(m,n)(z) - z \| \le \frac{\|z\|}{2}
$$
Now let $K_m$ be the compact subset of the basin of attraction
such that $F(m)(K_m)= B(r)$. Then we have that $G(m,n)^{-1}X_n
F(n)(K_m) \supset B(r/2)$ for any $n
>m$. It follows from part (b) of Lemma \ref{upper} that for every
$R>0$ there exist an $N \in \NN$ such that for $m \ge N$ we have
that
$$
B(R) \subset G(m)^{-1}\left(B(r/2)\right).
$$
Therefore $B(R) \subset \Psi(K_m)$ for large enough $m$, and thus
we have that $\Psi(\Omega)$ contains balls around the origin with
arbitrarily large radii. This completes the proof.

\end{document}